\documentclass{article}
\def\hS{ \hat{S} }
\def\Tel{ Tel^{\infty}}

\def\bA{ {\bf A} }
\def\bC{ {\bf C } }

\def\A{{ {\cal A} }}

\def\Spec{{ \mbox{Spec} }}

\def\1ox{{ \Omega^1_{\scriptstyle{X}} }}
\def\2ox{{ \Omega^2_{\scriptstyle{X}} }}

\def\ok1{{ \Omega^1_K }}
\def\ok2{{ \Omega^2_K }}
\def\Om{{ \Omega }}
\def\om{{ \omega  }}
\def\O{{ {\mathcal O} }}

\def\ra{{ \rightarrow }}
\def\da{{ \downarrow }}
\def\la{{ \leftarrow }}

\def\hra{{ \hookrightarrow }}
\def\da{{ \downarrow }}

\def\N{{ {\bf N} }}
\def\Q{{ {\bf Q} }}

\def\C{{ {\bf C} }}

\def\8{{ {\infty } }}

\def\s{{ {\sigma } }}
\def\G{{ \Gamma }}

\def\^{{ ^{\wedge} }}

\def\Z{{ {\bf Z } }}

\newtheorem{thm}{Theorem}
\newtheorem{cor}{Corollary}
\newtheorem{lem}{Lemma}
\newtheorem{prop}{Proposition}

\usepackage{amsmath}
\usepackage{amsfonts}
\usepackage{amscd}
\usepackage{amssymb}

\def\invlim{\varprojlim}
\def\dirlim{\varinjlim}

\def\hX{{ \hat{ X} }}

\def\H{{ \underline {H} }}
\def\C{{ {\cal C} }}

\usepackage{amsmath}
\usepackage{amsfonts}
\usepackage{amscd}
\usepackage{amssymb}
\def\invlim{\varprojlim}
\def\dirlim{\varinjlim}
\title{The Hyodo-Kato theorem for rational homotopy types}
\author{Minhyong Kim and Richard M. Hain}
\begin{document}
\maketitle
\section{Introduction}
Let $A$ be a complete discrete valuation ring with
perfect residue field $k$ of characteristic $p>0$
 and fraction field $F$. Denote by $W$  the ring
of Witt vectors of $k$ and by $K$ the fraction field
of $W$. Endow $A$ with the log structure
$A-{0} \hra A$ and
let $X$ be a  proper smooth connected fine saturated log scheme over $A$ with
 generic fiber $X^*$ and 
special fiber $Y$ which we assume to be of Cartier type
(\cite{HK} 2.12). The Hyodo-Kato theorem (\cite{HK} Theorem 5.1) says
$$H^i_{DR}(X^*)\simeq H^i_{cr}(Y/W)\otimes_W F$$
where the cohomology groups appearing in
the statement are  algebraic De Rham cohomology
on the left
and the crystalline cohomology of $Y$
with respect to $W$ on the right.

We will prove a version of this theorem in the
context of the unipotent rational homotopy types
defined in \cite{KH}. The definitions
 will be reviewed in the next section, but let us state the main result
here. If
$A_{cr}(Y)$ denotes the crystalline rational homotopy
type of $Y$ and $A_{DR}(X^*)$ the De Rham rational homotopy type
of $X^*$, then
\begin{thm}
$$A_{DR}(X^*)\simeq A_{cr}(Y)\otimes F$$
in the homotopy category of commutative differential
graded algebras (CDGA's) over $F$.
\end{thm}

In \cite{KH}, we proved this result essentially
when $A$ has the trivial
log structure.
Let $x$ be a point of $X$, $x^*$ the correponding point on
the generic fiber, and $y$ the reduction of $x$ to the special
fiber. It is simple to check that the augmentations
induced by $x^*$ and $y$ on the homotopy types are
compatible. Hence, we get

\begin{cor}
$$\pi_1^{dr}(X^*,x^*)\simeq \pi^{cr}(Y,y)\otimes_K F$$
as pro-unipotent algebraic groups over $F$.
\end{cor}

We also get results of the Artin-Mazur type \cite{AM} on
higher homotopy groups of simply connected varieties.
For this, let $X^*_1$ and $X^*_2$ be proper
 smooth geometrically-connected varieties over
a number field $E$ equipped with normal crossing divisors
$D^*_1$ and $D^*_2$.
Denote by $\O_{E,v}$  the localization at a prime $v$
of the ring of integers $\O_E$ of $E$.
Assume that $X_i^*$  extends to a proper flat scheme
$X_i$ over
$\O_{E,v}$ with special fibers $Y_i$ and that $D_i^*$
extends to a divisor $D_i$ on $X_i$
that is  relatively of normal crossing.
Further
assume that $D_i+Y_i$ is a strict normal crossing divisor.
 View
$X_i$ as a log scheme with the log structures
given by $D_i+Y_i$. $Y_i$ also has the induces log structure.
\begin{cor}
Suppose $Y_1 \simeq Y_2$ as log schemes.
Then for every embedding $\s:E \hra \bC$ such that
the $(X^*_i-D^*_i)\otimes_E \bC$
are simply-connected, their  higher rational homotopy
groups are isomorphic.
\end{cor}
The hypothesis of the corollary is satisfied for
example if
 $$X_1\otimes \O_{E,v}/m_v^2 \simeq X_2\otimes \O_{E,v}/m_v^2$$
by a map that preserves the divisors $D_i+Y_i$.

The proof of the theorem consists of choosing
explicit CDGA representatives for the homotopy types and a few
other intermediate complexes 
using embedding systems for crystalline cohomology,
constructing explicit multiplicative maps between them
and, thereby, `exorcising the derived category'
from the proof of the usual Hyodo-Kato isomorphism.

It is perhaps useful to think about the theorem in
the following  general context: Let $\A$ be the homotopy category of
non-negatively graded
CDGA's over some field and consider the forgetful
functor $$R:\A \ra \C$$ to the derived category of
complexes. Following the ideas of rational
homotopy theory, one should think of objects
$A$ in $\A$ as being rationally nilpotent topological 
spaces and $R(A)$ as being  like the homotopy
groups of the space. Given two objects
$A$, $B$ and an isomorphism
$g: R(A) \simeq R(B)$, we can then ask if there
is in fact an isomorphism $f:A \ra B$
such that $g=R(f)$. That is, we wish to prove a
Whitehead-type theorem showing that an isomorphism
of homotopy groups is induced by an actual isomorphism
of spaces. In our situation,
the usual Hyodo-Kato isomorphism says
that $R(A_{DR}(X^*)) \simeq R( A_{cr}(Y)\otimes_KF)$
and our theorem says that this
map can be lifted to a map of `spaces.'

In the process of making maps explicit, one encounters the
problem of assembling a collection of
maps in (an increasing sequence of) finite characteristics
into a single map of algebras. This is because the
Hyodo-Kato isomorphism depends on the fact
that repeated Frobenius twists of complexes computing
crystalline cohomology w.r.t. two different log structures
become increasingly close, i.e., isomorphic modulo
increasingly high powers of $p$. That is, an inverse system
of complexes related by the relative Frobenius
and reduction mediates the Hyodo-Kato isomorphism.
To deal with this
difficulty we use a differential graded algebra
version of homotopy direct limits, namely, the
mapping telescope of an inverse system of
maps. In fact, it will be necessary to,
employ the notion of a `twisted' inverse limit
introduced by Ogus to define 
`twisted' mapping telescopes and the key
definition of an `infinitely twisted'
mapping telescope (obtained by taking a direct limit
of twisted mapping telescopes).
Once one has this machinery in place,
the  isogenies that occur in the usual
proof fall into place nicely
as isomorphisms of $\infty$-twisted mapping
telescopes and the rest of the argument becomes  quite short.

We remark that the proof of theorem 1 appears to clarify  
the cohomological Hyodo-Kato isomorphism as well:
The rational homotopy types considered as complexes
compute the usual cohomology groups (crystalline
and De Rham). However, the proof given here,
while just a modification of those of
Hyodo-Kato and Ogus,
possesses an advantage over them in that the
maps involved are made completely explicit.
Furthermore, the proof makes clear that the
essence of the argument is local, although we
do not present it in an overtly local fashion.

One issue we do not deal with in this paper is
the dependence of the isomorphism stated in the
theorem on the choice of a uniformizer for $A$.
Although it is more or less clear
that one gets the same kind of
dependence as in the cohomological theorem,
we prefer to discuss this in a subsequent
paper together with a more detailed
study of  the monodromy operator and
period isomorphisms.

A few words about our convention: For the most part we
leave the log structure implicit and do not introduce
separate notation to indicate
their presence. The important
exception of course is in the proof of theorem 1.
A {\em multiplicative quasi-isomorphism}
 of CDGA's is just a quasi-isomorphism at the level
of complexes that respects the multiplicative
structure. A {\em quasi-equivalence}, on the other hand,
is an isomorphism in the homotopy category of CDGA's.

\section{Brief review}
For precise definitions, we refer the reader to
\cite{KH} sections 3 and 4, and the references therein.
In this section, we will just recall at the superficial
level some basic notation,
the objects that we will need, and their main properties.

For concreteness, we concentrate on
the situation described in the introduction.
Therefore, $k$ is endowed with the log structure
of the punctured point which is associated to
the pre-log structure $\N \ra k$ that sends
1 to 0.
 $Y$, then, is a smooth proper  fine saturated log scheme of Cartier type
over $k$ and $y:\Spec(k) \ra Y$ is a point of $Y$.
(We remind the reader that this means in particular
that it is a map of log schemes.)

The {\em unipotent crystalline rational homotopy type}
 $A_{cr}(Y)$ of $Y$ is defined by the following
formula:
$$A_{cr}(Y)=TW(W\om_y):=s_{TW}((\invlim \G (G (W\om_Y)))\otimes_W K)$$
The notation
is that given a pro-sheaf $L$,
$G(L)$ is a cosimplicial Godement resolution \cite{Go} for
the \'etale topology,
$\invlim$ goes from the category of inverse
systems $(M_n)_{n\in \N}$, where
$M_n$ is a  (cosimplicial) $W_n$-algebra, to the category of  complete
cosimplicial 
$W$-algebras, and finally,
$s_{TW}$ is Navarro-Aznar's `simple Thom-Whitney algebra'
functor \cite{NA}.   $W\om_Y$ is the pro-sheaf of CDGA's consisting
of the De Rham-Witt differential forms of $Y$. 
It is probably best at this point not to worry about
the precise construction and just remember that
$$TW(\cdot)=s_{TW}((\invlim \G (G (\cdot))\otimes)_W K)$$
is a functor from the isogeny category of pro-sheaves of
CDGA's over $W$ (that is, systems $(L_n)$
where a given level $L_n$ is a sheaf of CDGA's over $W_n$) 
to the homotopy category of
CDGA's over $K$. As motivation for the language of
homotopy types, we refer the reader to \cite{Mo}.

A basic proprty of
the functor $s_{TW}$ is that for any cosimplicial
CDGA $C$, $s_{TW}(C) \simeq s(C)$ in the derived category
of complexes,
where $s(\cdot)$ denotes the usual `associated
simple object' functor.
Ostensibly, the definition depends on the choice of a surjective
system of geometric
points on $Y$. 
Denote for a moment $A_{S}$ the rational
homotopy type constructed from a specific
system $S$ of geometric points.
Given two different systems $S_1$
and $S_2$,
we can find a third system $S_3$ dominating both
$$\begin{array}{ccccc}
& & S_3 & & \\
& \swarrow & & \searrow \\
S_1& & & & S_2
\end{array}$$
which therefore gives us multiplicative
quasi-isomorphisms
$$\begin{array}{ccccc}
& & A_{S_3} & & \\
& \nearrow & & \nwarrow \\
A_{S_1}& & & & A_{S_2}
\end{array}$$
Hence, the isomorphism class of $A_{cr}(Y)$ in
the homotopy category of CDGA's
is well-defined. Furthermore,
if $S_4$ were chosen instead of $S_3$, then we could
dominate these two by a third $S_5$ fitting into a commutative
diagram
$$\begin{array}{ccccc}
& & S_3 & & \\
& \swarrow &\uparrow & \searrow \\
S_1& \la& S_5& \ra& S_2\\
& \nwarrow & \da & \nearrow \\
& & S_4 & & 
\end{array}$$
leading us to a commutative diagram
$$\begin{array}{ccccc}
& & A_{S_3} & & \\
& \nearrow &\downarrow & \nwarrow \\
A_{S_1}& \ra& A_{S_5}& \la& A_{S_2}\\
& \searrow & \uparrow & \swarrow \\
& & A_{S_4} & & 
\end{array}$$
Therefore, the isomorphism $A_{S_1} \ra A_{S_2}$
is canonical 
in the homotopy category. This justifies omitting the
system of points from the notation for the homotopy
type. We will omit similar obvious arguments in
a few other places of the papers.

The (unipotent) De Rham rational homotopy type of $X^*$
is defined by
$$A_{DR}(X)=TW(\Om_{X^*/F} ):=s_{TW}\G (G (\Om_{X^*/F}))$$
where the Godement resolution is now taken on $X^*$.
But GAGA implies that
$$A_{DR}(X)\simeq s_{TW}((\invlim \G (G (\Om_{\hX/A})))\otimes_A F)$$
where  the formal completion $\Om_{\hX/A} $ of the De Rham complex
of $X/A$ is regarded as a pro-sheaf on $Y$.
It is the latter object that we will compare to 
$A_{cr}$.

 $A_{cr}$ can be computed using
the `crystalline complex' associated to an
embedding system (\cite{HK} 2.18): Let $(Y_.,Z_.)$ be a pair of simplicial
schemes that fit into a diagram
$$\begin{array}{ccc}
Y_. & \hra & Z_. \\
\da & & \\
Y& & 
\end{array}
$$
where $Y_.$ is a simplicial hypercovering (that
satisfies cohomological descent for the \'etale topology)  which we equip
with the log structure pulled back from that of $Y$, 
$Y_.\hra Z_.$ is a closed embedding of formal log schemes,
and $Z_.$ is smooth formal log scheme over $W$. Then the associated
crystalline complex is by definition
$$C(Y_.,Z_.):=\Om_{Z_./W} \otimes_{\O_{Z_.}} D_{Y_.}(Z_.)$$
Here, $D_{Y_.}(Z_.)$ is the divided power envelop of $Y_.$
in $Z_.$. Regard $C(Y_.,Z_.)$ as a simplicial pro-CDGA on
$Y_.$. Then we also defined $TW(C(Y_.,Z_.))$ in this
setting and we have a quasi-equivalence (i.e.,
an isomorphism in the homotopy category):
$$TW(C(Y_.,Z_.))\simeq TW(W\om_Y)$$
The functor $TW$ (Thom-Whitney) can be defined in a more general
setting, as we have already done in the definition of
$A_{DR}$ for example, but also for crystalline
complexes over more general bases.
So if $Y/S$ is a fine saturated smooth log scheme over an affine
base $S=\Spec(R)$ and $S\hra T=\Spec(B)$
is an exact immersion of
formal log schemes where $B$ is of characteristic zero,
then to any embedding system $(Y.,Z_.)$ as above
with $Z_.$ smooth over $B$, we can associate the
crystalline complex $C(Y_.,Z_.)$ and the Thom-Whitney
algebra $TW(C(Y_.,Z_.))$ which ends up as a CDGA over $B\otimes \Q$.
This algebra is in fact canonically independent of
the embedding system: Any two embedding systems
can be dominated by a third giving rise to a quasi-equivalence
which, in turn, is independent in the homotopy
category of the dominating system.

\section{Mapping telescopes}

Now we will go on to define the {\em (twisted) mapping telescopes}
of inverse systems of complete $W$-algebras.

First, if $f: A\ra B$ is a map of (complete) CDGA's over $W$,
we define the {\em mapping cylinder} that fits into a diagram
$$\begin{array}{ccc}
 & & Cyl(f) \\
 & \nearrow & \da \\
A & \ra & B
\end{array}$$
as follows:
Denote by $I$ the formal divided power De Rham complex of $W<<x>>$, which
therefore is concentrated in degrees  0 and 1.
Consider the two augmentation maps
$e_0, e_p: I \ra \Z$ that evaluate functions in $I^0$
at 0 and $p$, respectively, and sends $I^1$ to zero.
Hence $B\otimes I$ is likewise equipped with two maps to $B$
which we also denote by $e_0$ and $e_p$. Here and henceforward,
all tensor products are topological. The mapping
cylinder $Cyl(f)$ of $f$ is defined to be the
subalgebra of $A\oplus B\otimes I$ consisting of elements
$(a,b)$ such that $f(a)=e_0(b)$. Notice that the
map $Cyl(f) \ra B$ induced by $e_p:B\otimes I \ra B$
becomes surjective after tensoring with
$\Q$: for any $b\in B$, $(0, xb/p)$ is in 
$Cyl(f)\otimes \Q$ and maps to $b$. On the other hand,
integrally, we can only say that if $b\in B$, then
$pb$ is in the image of $e_p$.
The map $A \ra Cyl(f)$ given by $a \mapsto (a, f(a))$
induces a q.i. In fact, the projection
$Cyl(f) \ra A$ to the first component is a chain homotopy
inverse.
Now, given an inverse system of maps 
$$A_.: \cdots \ra A_3 \ra A_2 \ra A_1$$
indexed by the positive integers, it is clear how to
construct the mapping telescope $Tel(A_.)$.
It is the inverse limit of the inverse system defined inductively by
putting $Tel_1=A_1$ and 
$$Tel_{i+1}:=Cyl(A_{i+1} \ra Tel_i)$$
where the map $A_{i+1} \ra Tel_i$ is
the composite $A_{i+1} \ra A_i \ra Tel_i$.
The construction of $Tel$ is clearly functorial
for inverse systems of CDGA's.
In particular, one can apply the construction to
inverse systems of cosimplicial CGDA's or
bi-cosimplicial CDGA's.
Also, there is  a map of inverse systems
$A_. \ra Tel_.$ giving rise to a 
functorial map 
$\invlim A_. \ra Tel$.

We will also need the `twisted' inverse limit
construction of Ogus (\cite{Og} p.203):
Given an inverse system $L_.$ of $W$-modules indexed by $\N$
and an integer $m$, let $\invlim^m L_.$ be the collection of
elements $(a_i)$, $a_i\in L_i$, such that
$a_i \mapsto p^m a_{i-1}$.
If $n\geq m$, then we have natural maps
$\invlim^m L_.\ra \invlim^n L_.$  given by
$(a_i) \mapsto (p^{(n-m)i}a_i)$. That is, the twisted
inverse limits form a directed system.

Using these, we can also construct the $m-$twisted
telescope $$Tel^m(A_.)=\invlim^m Tel_i(A_.)$$
as well as the infinitely twisted telescope
$$Tel^{\infty} (A_.)=\dirlim_m Tel^m (A_.)$$
Clearly, we have functorial maps
$\invlim^m A_. \ra  Tel^m (A_.)$
and
$$\dirlim_m \invlim^m A_. \ra  Tel^{\infty} (A_.)$$

\begin{lem}
Let $L_.\ra K_.$ be a map of inverse systems of  
((bi-)cosimplicial) CDGA's
which is a quasi-isomorphism at each level. That is,
each $L_n \ra K_n$ is a quasi-isomorphism.
Suppose, furthermore, that the transition maps for
both systems (denoted $\pi$) have the following property:
there exists some $m$ such that for each element
$x$ of level $i$, there exists an element $y$ of level $i+1$
such that $\pi(y)=p^mx$.
Then
$$\dirlim_n \invlim^n L_. \ra \dirlim_n \invlim^n K_.$$
is a quasi-isomorphim.
\end{lem}
{\em Proof.} It is straightforward to check that
the twisted inverse limit commutes with the cone construction,
so that, for each $n$, we have an exact triangle
$$0\ra \invlim^n L_. \ra  \invlim^n K_.\ra \invlim^n C_.$$
where $C_.$ is the cone of $L_. \ra K_.$. Since direct limits
commute with cohomology, we need only show that
$$\dirlim H(\invlim^n C_.)=0$$
In fact, the transition map
$$H(\invlim^n C_.) \ra H(\invlim^{n+m} C_.)$$
is zero.
To see this, let the cocycle $(c_i)$ represent an element of
$H(\invlim^n C_.)$. We will show that
$(p^{mi}c_i) \in \invlim^{n+m} C_.$
is a coboundary. Construct an element $(b_i) \in \invlim^{n+m} C_.$
such that $d(b_i)=(p^{mi}c_i)$ inductively as follows:
assume we have constructed up to $b_j$. That is, for $i\leq j$,
$db_i=p^{mi}c_i$ and $\pi (b_j)=p^{n+m}b_{j-1}$.
Since each $C_i$ is acyclic, there exists an
$x_{j+1} \in C_{j+1}$ such that $dx_{j+1}=c_{j+1}$.
Hence, $d\pi (p^{mj}x_{j+1})=p^{mj}\pi (c_{j+1})=p^{mj+n}c_j=p^ndb_j$
and $d (\pi (p^{mj}x_{j+1})-p^nb_j)=0$. Again by
acyclicity in level $j$, we can then find an $a_j$
such that $da_j=\pi (p^{mj}x_{j+1})-p^nb_j$. Multiplying
by $p^m$, we get
$$\pi (p^{(j+1)m}x_{j+1})-p^{n+m}b_j=d(p^ma_j)=d\pi (a_{j+1})$$
for some $a_{j+1}$. Now put
$$b_{j+1}:=
p^{(j+1)m}x_{j+1}-da_{j+1}$$
Then $\pi(b_{j+1})=p^{n+m}b_j$ and $db_{j+1}=p^{(j+1)m}c_{j+1}$,
so we are done.
\medskip

We emphasize that
$\Tel$ is a functor from inverse systems of
 ((bi-)cosimplicial) CDGA's to ((bi-)cosimplicial) CDGA's
and that the maps
$$\invlim A_.  \ra \dirlim_m \invlim^m A_. \ra  Tel^{\infty} (A_.)$$
are multiplicative.

\section{Proof of theorem}

In this section, it will be useful to employ the following
notation: Given a $p$-adic formal scheme $S$ and an integer
$m$, 
$S/m$ denotes $S\otimes_{\Z_p}\Z/p^{v_p(m)}$, where
$v_p$ is the $p$-adic valuation for which $v_p(p)=1$.

 We denote by
$S_1$ and $S_2$ the  scheme $\Spec W [t]$ equipped with
the log structures associated to the pre-log structures
$\N \ra W[t]$ that send 1 to $0$ and $t$, respectively.
So both restrict to the log structure of the punctured point on
$\Spec (k)$.
$W[t]$ carries the Frobenius $\s$ given by the usual Frobenius
on $W$ and sending $t$ to $t^p$. Denote by
$W<<t>>$, the $p$-adic completion of
the divided power polynomial algebra $W<t>$ over $W$
which therefore  also carries log structures induced by
those of $S_1$ and $S_2$. Denote these formal log schemes
by $\hS_1$ and $\hS_2$.  The Frobenius $\s$ of $W[t]$ 
naturally extends to $W<<t>>$. We will
denote by $S_i^{(n)}$ the 
scheme $\Spec (W[t])$ equipped with the log structure
of $S_i$ pulled back by the $n$-th iterate of the Frobenius.
Hence, we see that $S_1^{(n)}\simeq S_1$ while
the log structure on $S_2^{(n)}$ is associated to
the pre-log structure $\N \ra W[t]$ that sends
$1$ to $t^{p^n}$. We will also use the notation
$\hS_i^{(n)}$ in the obvious manner.

The difficult part of the Hyodo-Kato isomorphism says
that if $Y$ is a smooth fine proper log scheme of
Cartier type, then
$$H^i(Y/S_1)\otimes_{\Z} \Q \simeq H^i(Y/S_2)\otimes_{\Z} \Q$$
that is, the crystalline cohomology is
the same for the two log structures, up to isogeny.

We recall  the proof, taking care to make some maps
more explicit. 

We start out by choosing  embedding systems $Y_.\hra Z^1_.$ and
$Y_. \hra Z^2_.$ for $Y$ w.r.t. the two log schemes
$S_1 $ and $S_2$ that admit Frobenius lifts. This
notion requires a brief explanation: Denote by
$S_a^{(n)}$ the Frobenius twisted log
schemes introduced above. There are
then maps $S_a \ra S_a^{(n)}$ 
 induced from the monoid map
$\N \ra \N$, $1\mapsto p^n$. Now,
$(Z^a_.)^{(n)}$ is a smooth simplicial log scheme
over $S_a^{(n)}$, and a Frobenius lift refers to
a map $F: Z^a_.  \ra  (Z^a_.)^{(1)}$ that fits into the commutative diagram
$$\begin{array}{ccc}
Z^a_. & \stackrel{F}{\ra} & (Z^a_.)^{(1)} \\
\da & & \da \\
S_a & \ra & S_a^{(1)}\\
\end{array}$$
Such embedding systems can be constructed,
for example, as follows: $Y_0$ in both cases is just
the disjoint union of the elements of an affine open
covering of $Y$ and $Z^1_0$ and $Z^2_0$ are
smooth liftings of $Y_0$ to $S_1$ and $S_2$,
respectively. Then the $Y_i=Y_1\times_Y Y_1 \times \cdots \times Y_1$
($i+1$-times)  embed diagonally into
 $$Z^a_i=Z^a_0\times_{S_a} Z^a_0 \times_{S_a} \cdots \times_{S_a} Z^a_0$$
($i+1$-times) for $a=1,2$ and come together to form a simplicial
hypercovering. By the affine smoothness, we easily
get the Frobenius lifts  in cosimplicial level 0 and
then the other levels by taking products.

Let $C_1$ and $C_2$ be crystalline complexes for the
two different embedding systems, and form the
bi-cosimplicial algebras over $W<<t>>$
$$B_a:=\invlim \Gamma (G(C_a))$$
for $a=1,2$.
Denote by $B^{(n)}_a$ the pull-back of $B_a$ via
$\s^n$, which therefore arise from crystalline complexes
for $(Z^a_.)^{(n)}$, the $n$-th Frobenius pull-back of
$Z^a_.$, w.r.t. $S_a^{(n)}$.
 In particular, when we adjoin divided powers and reduce
mod $p^n!$, the two log structures are isomorphic,
so we get  isomorphisms in the derived category of
bi-cosimplicial
complexes
$$B^{(n)}_1/(p^n)! \simeq B^{(n)}_2/(p^n)!$$
We can get these maps to  fit into commutative diagrams
of bi-cosimplicial CDGA's
$$\begin{array}{ccccc}
B^{(n+1)}_1/(p^n)! &\ra & L_{n+1}&\la &B^{(n+1)}_2/(p^n)! \\
\da \scriptstyle{\phi} &&\da & & \da \scriptstyle{\phi}\\
B^{(n)}_1 /(p^n)! &\ra &L_n &\la &B^{(n)}_2/(p^n)!
\end{array} $$
where the horizontal arrows are multiplicative quasi-isomorphisms
 of CDGA's (as opposed to maps in the derived category of
complexes) constructed as follows. The vertical arrows
on either end are induced by the Frobenius lifts.
$L_n$ is the crystalline complex for $Y^{(n)}$ relative to
$\hS_1^{(n)}/p^n!=\hS_2^{(n)}/p^n!$ 
computed using the `diagonal' embedding system
$$Y^{(n)}_.\hra (Z^1_.)^{(n)} \times_{\hS^{(n)}_1/p^n!} (Z^2_.)^{(n)} $$ 
The Frobenius lifts for $Z^1_.$ and $Z^2_.$ determine one
for each product giving the maps $L_{n+1}\ra L_n$.
The horizonal maps are then induced by the projections.
We stress that all the maps are given by
pull-backs of differential forms, and hence, are multiplicative.

On the other hand, the iterates of
the Frobenius induce  multiplicative maps of inverse  systems
$$B_a^{(.)}/p^.! \ra B_a/p^.!$$ where the
target is just the inverse system given by
the reductions of $B_a$ with the natural projections connecting them.

Lemma 1 from the previous section admits the following
corollaries:

\begin{cor}
The following arrows are quasi-isomorphisms:
$$Tel^{\infty} (B^{(.)}_1) \ra \Tel (L_.) \la \Tel (B^{(.)}_2) $$
\end{cor}
\begin{cor}
$$\dirlim_m \invlim^m (B_a/p^.!) \ra \Tel (B_a/p^.!) $$
is a quasi-isomorphism.
\end{cor}
On the other hand, we have the 
\begin{prop}
The map
$$B^{(.)}_a/p^.! \ra B_a/p^.!$$
induces a quasi-isomorphism
$$\Tel (B^{(.)}_a/p^.! ) \simeq \Tel (B_a/p^.! )$$
for $a=1,2$. 
\end{prop}
{\em Proof.}
We will omit the subscript $a$ from the notation.
Hence, $S$ also refers to either $S_1$ or $S_2$.

If $C_n$ is the cone of $B^{(n)}/p^n! \ra B/p^n!$
then $Tel_n (C_.)$ is the cone of
$Tel_n(B^{(.)}/p^.! ) \ra  Tel_n (B/p^.!)$
That is, we have an exact triangle:
$$0\ra Tel_n(B^{(.)}/p^. ) \ra  Tel_n (B/p^.!) \ra Tel_n (C_.)$$

The key point is the following
\begin{lem}

$H^i(C_n)=H^i(Tel_n(C_.))$
is killed by $p^{2in}$ for $i>0$
and $p^n$ for $i=0$.
\end{lem}
{\em Proof of lemma.}
Since the algebras are constructed out of the
stalks of Godements resolutions, the statement is
local on $Y$. So we may assume that the embedding system
is just a smooth $S$ lift $X$ and that we also
have a Frobenius lift $f$. Thus, $f$ induces pull-back maps
$$\phi^n:\Om_{X^{(n)}/S^{(n)}} \ra \Om_{X/S}$$
It suffices to show that the cone $C_n$  of this map
mod $p^n!$
has $i$-th cohomology killed by $p^{2in}$
for $i>0$ and $p^n$ for $i=0$.
We just give the argument for $i>0$ since
the $i=0$ case is an obvious modification.

Recall the map $F$ given by $\phi/p^i$ in degree
$i$. The definition works just as in the usual
case (\cite{Il} 0.2.3.3)  using the $W$-flatness of the sheaf of differentials.
 $F$ induces an injection (\cite{IR} III.1.5.4):
$$F:\Om^i_{X^{(1)}/S^{(1)}}/[p^n(\Om^i_{X^{(1)}/S^{(1)}})+pd(\Om^{i-1}_{X^{(1)}/S^{(1)}})]
\hra \Om^i_{X/S}/[p^n(\Om^i_{X/S})+d(\Om^{i-1}_{X/S})]$$
This  is an immediate consequence of the Cartier isomorphism,
again as in the classical case.

Note the following corollary: given $x \in \Om^i_{X^{(1)}/S^{(1)}}$,
if $F(x)=p^my$ for some $y$, then $x=p^mz$ for
some $z$. This is obvious for $m=1$ from the above injection.
Assume it for $m-1$. $F(x)=p^my$ in any case implies
$x=p^mw+pdu$. But then
$F(pdu)=F(x)-F(p^mw)=p^m(y-F(w))$ so $F(du)$ is divisible
by $p^{m-1}$. Therefore, $du=p^{m-1}v$ and $x$ is divisible
by $p^m$.

By iterating the argument, we also see that if
$x\in \Om^i_{X^{(n)}/S^{(n)}}$ and
$F^n(x)$ is divisible by $p^m$, then so is $x$.

Now, let 
$$(a,b)\in\Om^i_{X/S}\oplus \Om^{i+1}_{X^{(n)}/S^{(n)}}$$ 
represent a cocycle in
$$C^i=(\Om^i_{X/S}\oplus \Om^{i+1}_{X^{(n)}/S^{(n)}})\otimes \Z/p^n!$$
Then
$\phi^n(b)-da=p^n!x$ and $db=p^n!y$.
Write the first equality as
$$F^n(p^{in}b)=p^n!x+da$$
and apply the injection above to get
$$p^{in}b=p^n!w+pds=p^n!w+dc$$
for some $w$ and $c$. We also get
$$d(\phi^n(c)-p^{in}a)=p^{in}\phi^n(b)-p^n!p^{in}F^n(w)-p^{in}da=
p^n!p^{in}(x-F^n(w))$$
That is,
$$\phi^n(c)-p^{in}a \in d^{-1}(p^n!p^{in}\Om^{i+1}_{X/S})$$
Recall the  formula
$$ d^{-1}(p^n\Om_{X/S}^{i+1})=
\Sigma_{0\leq k\leq n}p^kF^{n-k}(\Om_{X^{(n-k)}/S^{(n-k)}}^i)+
\sum_{0\leq k\leq n-1}F^k(d\Om_{X^{(k)}/S^{(k)}}^{i-1})$$
whose proof also follows \cite{Il} 0.2.3.13 verbatim, as pointed out
by Jannsen
in the appendix of \cite{Hyo1}. Thus, we have
$$\begin{array}{c}
\phi^n(c)-p^{in}a=
p^n!p^{in}z_0+p^n!p^{-1}p^{in}F(z_1)+
\cdots+p^n!p^{-n}p^{in}F^n(z_n)+\cdots\\
 F^{N}(z_N)+du_0+F(du_1)
+\cdots+F^{N-1}(du_{N-1})\end{array}$$
where $N=v_p(p^n!)+in$. We examine the terms in this
equality after multiplying by another $p^{in}$. We get that
$\phi^n(p^{in}c)-p^{2in}a$  is of the form
$p^n!p^{in}z+\phi^n(l)+du$. If we apply the differential
$d$, we get that $d \phi^n(l)=\phi^n(dl)$ is divisible
by $p^{in}p^n!$. Thus, $F^n(dl)$ is divisible by $p^n!$,
and hence, so is $dl$. So we get
$d(p^{in}c-l)=p^{2in}b $ (mod $p^n!$)
and $\phi^n(p^{in}c-l)-du=p^{2in}a $ (mod $p^n!$).
Hence, the class of $(p^{2in}a, p^{2in}b)$ (mod $p^n!$)
is a coboundary. This proves the lemma.
\bigskip

For all our considerations, we may assume that the dimension
$d$ of $Y$ is positive. Therefore, we see that all the
cohomology of $C_n$ is killed by $p^{2dn}$.
The remainder of the proof of the proposition
is as in the proof of
Lemma 1 : Let $(c_n)$ be a cocycle in $Tel^m (C_.)$. Then
$(p^{(2d+1)n}c_n) \in Tel^{(m+2d+1)}(C_.)$ is a coboundary.

\medskip

{\bf Remark:} The fact that the relative Frobenius
is an isogeny has been the source of
many important theorems on crystalline cohomology
starting with the theorem of Berthelot and Ogus.
Here, it is rather mysteriously manifested
in the isomorphism of infinitely twisted mapping
telescopes.
\medskip

Thus we have arrived at an explicit isomorphism in
the derived category from
$$\dirlim^m \invlim^mB_1/p^.!$$ to $$\dirlim^m \invlim^m B_2/p^.!$$
mediated by the following system of arrows all of
which are quasi-isomorphisms:
$$\begin{array}{ccccccccc}   
& & & &\Tel (L_.)& & & &  \\
& & & \nearrow & & \nwarrow & & &\\
& &\Tel (B^{(.)}_1/p^.!)& & & & \Tel (B^{(.)}_2/p^.!) & &\\
& & \da & & & & \da & & \\
 \dirlim^m \invlim^mB_1/p^.!&\ra&\Tel(B_1/p^.!)& & & & \Tel(B_2/p^.!)& \la &\dirlim^m \invlim^mB_2/p^.!
\end{array}$$
On the other hand, as has already been
pointed out, it is trivial to check
that if $L_.$ is an inverse system of $W$-algebras,
then $\dirlim^m \invlim^m L_. $ is a $W$ algebra,
even though each $\invlim^m L_.$ separately is not.
Also, maps in the category of inverse systems of
algebras
induce algebra homomorphisms for this double limit.
Therefore, we see that the isomorphism constructed
above is actually an isomorphism in the homotopy
category of algebras.

Now consider the natural map
$$\invlim B_a/p^.! \ra \dirlim^m \invlim^m B_a/p^.!$$
By Ogus (\cite{Og} Lemma 18), we have an isomorphism
$$(\invlim B_a/p^.! )\otimes \Q \simeq 
(\dirlim^m\invlim^m B_a/p^.!)\otimes \Q$$
Thus, we have constructed a quasi-equivalence from
$$(\invlim B_1/p^i!) \otimes \Q$$
to
$$(\invlim B_2/p^i!) \otimes \Q$$ 
Applying the Thom-Whitney simple object functor
twice, we get a quasi-equivalence
$$TW( C(Y_., Z_.^1) \simeq TW( C(Y_.,Z^2_.)) $$
It is straightforward to check that
this map is independent (in the homotopy category)
 of the various choices
made.

The rest of the argument is as in \cite{KH}, section 7,
and amounts to a `Berthelot-Ogus' type argument \cite{BO}.
Choose a uniformizer $\pi$ of $A$ which therefore
determines a presentation $A\simeq W[t]/(f(t))$,
where $f(t)$ is an Eisentein polynomial of degree $e=[F:K]$.
Let $R$ be the $p$-adic completion of the
divided power envelop of $(f(t),p)$ inside $W[t]$. 
Thus $R$ is also the completed DP envelop
of the ideal $(t^e)$.
We have a natural map $g:R\ra W<<t>>$. On the other hand,
if   $r$ is such that $p^r\geq e$,
 then the map
$\s^r:W<<t>>\ra W<<t>>$  factors through
$W<<t>> \ra R \ra W<<t>>$. Composing in the
other direction
$R \ra W<<t>> \ra R$ is, by definition, the
Frobenius map of $R$. 
Let $Y'=X\otimes A/p$ with the induced log structure.
Also, give $\Spec(R)$ the log structure induces from $S_2$.
We have the commutative diagram
$$\begin{array}{ccccc}
Y'& \ra &Y & \hra& Y' \\
\da & & \da & & \da \\
\Spec(A/p) & \ra & \Spec(k) & \hra & \Spec(A/p)\\
\da &  &\da & & \da \\
\Spec(R) & \ra & \hS_2 & \ra & \Spec(R)
\end{array}$$
where the composite of the horizontal arrows are all
the $r$-th iterate of the Frobenius.
We  choose crystalline complexes for $Y$ and $Y'$
as
follows.
Construct first embedding systems for $Y$ and $Y'$
w.r.t. $S_2$ that fit into a diagram
$$\begin{array}{ccccc}
Y_. & \hra & Y'_. & \hra &Z_. \\
\da & & \da & &  \\
Y&\hra & Y' &  & \\
\end{array} $$
where $Z_.$ is smooth over $S_2$  and the
left hand square is cartesian. We can also arrange for
$Z_.$ to admit a Frobenius lift compatible with the
Frobenius of $W[t]$. Then
$C=\Om_{Z_./W[t]}\otimes W<<t>>$ and
$C'=\Om_{Z_./W[t]}\otimes R$ are crystalline complexes
for $Y$ and $Y'$ and we can regard both
as simplicial sheaves of CDGA's on $Y_.$.

Denote by $C^{(r)}$ (resp. $(C')^{(r)}$ ) the pull-back of $C$
(resp. $C'$)
by the $r$-th power of the Frobenius map of $W<<t>>$
(resp. $R$).
Then the big diagram above implies that
$$(C')^{(r)}\simeq C^{(r)}\otimes_{W<<t>>} R$$
in the homotopy category of
sheaves of CDGA's on $Y_.$. On the other hand,
the Frobenius lifts induce maps
$C^{(r)} \ra C$ and $(C')^{(r)} \ra C'$.
So we have  maps  of sheaves of
CDGA's
$$C\otimes R \leftarrow C^{(r)}\otimes R \simeq (C')^{(r)} \ra C'$$
and taking Thom-Whitney algebras, we have multiplicative
maps
$$TW(C)\otimes (R\otimes \Q)  \leftarrow TW(C^{(r)})\otimes (R\otimes \Q)
\simeq TW((C')^{(r)}) \ra TW(C')$$
of CDGA's over $R\otimes \Q$.
From the fact that the relative Frobenius is an isogeny 
(\cite{HK} Proposition 2.24), we
know that all these maps are quasi-equivalences.
Now we tensor with the quotient map
$R\otimes \Q \ra A\otimes \Q=F$ to get
$$TW(C) \otimes_{W<<t>>\otimes \Q} F \simeq TW(C')\otimes_{R\otimes Q} F$$
By using the fact that
$C'\otimes A$ is the crystalline complex
associated to an embedding system for $Y'$
w.r.t. $A$ which is also true of
$\Om_{\hX/A}$, we get
$$TW(C')\otimes F \simeq TW (\Om_{\hX/A}) \simeq TW (\Om_{X^*/F})$$
On the other hand, $C$ is quasi-equivalent to $C(Y.,Z^2_.)$
from the previous section, so that
$TW(C)\simeq TW(C(Y.,Z^2_.))\simeq TW(C(Y.,Z^1_.))$
and $C(Y.,Z^1_.)$ is quasi-equivalent to
the base change to $W<<t>>$ of a crystalline complex
for $Y$ w.r.t. $W$ so we get
$$TW(C)\otimes F \simeq TW(W\om_Y)\otimes_K F$$
giving us the desired quasi-equivalence
$$TW(W\om_Y)\otimes_K F \simeq TW (\Om_{X*/F})$$
This is the isomorphism of homotopy types stated
in the theorem.
\medskip

{\em Proof of corollary 1}

We need to discuss basepoints.
We start with a careful discussion of the
basepoint for $W\om_Y$. 
If $y:\Spec(k) \ra Y$ is a point of $Y$, then there is a map
$W\om_Y \ra y_*(W)$. It is given by
 $0$ in positive degrees
and the canonical map $e_y:W\O_Y \ra W(k)$ induced by
$\O_Y \ra k$. This map induces
$TW(W\om_Y) \ra TW (y_*(W))$. However, a simple examination of
the definition yields the following description
of the degree zero term $TW^0 (y_*(W))$: It consists of collections
$(f_n)$ where $f_n$ is a function on $\bA^n_K$
with the property that $\partial_i(f_{n+1})=f_n$ for any
$i$.  The map $(f_n) \ra f_0$
yields a quasi-isomorphism
$TW (y_*(W)) \simeq K$.
We obtain thereby the augmentation map
$TW(W\om_Y) \ra K$. Let's describe this map explicitly.
The degree zero term of $TW(W\om_Y)$ consists of
compatible sequences
(\cite{KH} section 3) $(a_n)$, 
$a_n \in \prod_{\bar{p}}(W\O_Y)_{\bar{p}}\otimes \O_{\bA^n_K}$
where $\bar{p}$ is an $n+1$-tuple of points in
$Y$. Thus, $a_0$ is just an element of $\prod_{p}(W\O_Y)_{p}$,
and the map $TW(W\om_Y) \ra K$ is zero in positive degrees
while in degree zero it sends $(a_n)$ to $e_y(a_0)$.
Now, let $(Y_.,Z_.)$
be an embedding system for which
$Y_0$ is the disjoint union of an affine open cover of $Y$
and $Z_0$ is a smooth lifting of $Y_0$ which admits a
Frobenius lift $F$. By the construction at the beginning of
the section, for example, such embedding systems exist.
The base point $y \in Y$ lifts to $ Y_0$ and a $W$ point
$z$ of $Z=Z_0$. Locally, we can express
$W_n\O_Y$ as the cohomology $\H^0(\Om_{Z_n/W_n})$
and then we get the map
$\O_Z \ra W\O_Y$ that sends $a$ to the sequence
$(F^n(a) (\mbox{mod} p^n))_n$ which is the degree zero component of the
quasi-isomorphism
$\Om_{Z/W} \ra W\om_Y$. Let $\Om_{Z/W}\otimes D_y(Z)$ be the
divided power envelop of $y$ in $Z$. Then we have
$(\H^0( \Om_{Z_n/W_n}\otimes D_y(Z_n)))_n \simeq W(k)$. The augmentation
$e_z:\O_Z \ra W$ given by evaluation at $z$
then fits into a commutative diagram
$$\begin{array}{ccccc}
\O_Z &  \ra &(\H^0(\Om_{Z_n/W_n}))_n & \simeq & W\O_Y\\
\da & & \da & &\da\\
W & \ra & (\H^0( \Om_{Z_n/W_n}\otimes D_y(Z_n)))_n & \simeq & W
\end{array}$$ where the first arrow in the bottom row
sends $a \in W$ to the inverse system $(\s^n(a) (\mbox{mod} p^n))_n$
This implies that
$$(TW(C(Y_.,Z_.)),e_z) \simeq (TW(W\om_Y), e_y)$$
as augmented algebras.

Now, assume  we are in the situation at the end of
section 2 where
$V \ra \Spec(R)$ is a smooth proper connected fine saturated
log scheme over an affine base and $\Spec(R) \hra \Spec(B)$
is an exact topologically nilpotent
immersion where $B$ is of characteristic zero.
Assume also that we are given a point $v:\Spec(R) \ra V$.
Then we can always find an embedding system
$(V_.,Z_.)$ with the property that $v$ lifts to $V_0$
and to a point $z: B \ra Z_0$. This point allows us to
put an augmentation
$TW(C( V_.,Z_.)) \ra B \otimes \Q$ on the
Thom-Whitney algebra and a product construction shows that
for any two choices of embedding systems, there
is a homotopy equivalence between the Thom-Whitney
algebras which is compatible with the augmentation.
Similarly for different choices of liftings of
the point $y$.
From this, applied to the various embedding systems
that occur in the proof of the theorem,
 one easily deduces that the homotopy
equivalence
$$TW(C(Y_.,Z_.)) \otimes F \simeq TW(\Om_{X^*})$$
takes the augmentation induced by $e_z$ to that induced by
evaluation at $x^*$.

\medskip

{\em Proof of Corollary 2.}

The isomorphism classes of the higher rational homotopy groups are
determined by their dimension, and this dimension can be
computed in any complex embedding of $E$ or after base
change to the completion $E_v$ of $E$ w.r.t. $v$. The assumptions
imply that the special fibers $Y_1$ and $Y_2$
 are isomorphic smooth log schemes.
Thus, $TW (\om_{Y_2}) \simeq TW(\om_{Y_2})$, which implies the
quasi-equivalence of $TW(\Om_{X_1^*}(\log D_1^*))\otimes E_v$
and $TW(\Om_{X_2^*}(\log D_2^*))\otimes E_v$. 
Thus, their bar complexes 
are quasi-equivalent, giving isomorphisms of their cohomology groups,
i.e., the higher De Rham homotopy groups of $X_1^*-D_1^*$ and $
X_2^*-D_2^*$ \cite{Wo}.
\bigskip

{\bf Acknowledgement:} Both authors were supported in part
by grants from the National Science Foundation.

{\footnotesize
M.K.: DEPARTMENT OF MATHEMATICS, UNIVERSITY OF ARIZONA, TUCSON, AZ 85721}

{\footnotesize EMAIL: kim@math.arizona.edu}

{\footnotesize R.H.: DEPARTMENT OF MATHEMATICS, DUKE UNIVERSITY, DURHAM, NC 27708} 

{\footnotesize  EMAIL: hain@math.duke.edu}
 

\begin{thebibliography}{30}
\bibitem{AM}
Artin, M.; Mazur, B. Etale homotopy. Reprint of the 1969 original. Lecture
Notes in Mathematics, 100. Springer-Verlag, Berlin, 1986. iv+169 pp. 

\bibitem{BO} Berthelot, P.; Ogus, A. $F$-isocrystals and de Rham cohomology. I. Invent. Math. 72 (1983), no. 2, 159--199. 
\bibitem{Go}
 Godement, Roger Topologie al\'ebrique et th\'eorie des faisceaux.
(French) Actualit\'es Sci. Ind. No. 1252. Publ. Math. Univ. Strasbourg. No. 13 Hermann,
Paris 1958
\bibitem{Hyo1}
Hyodo, Osamu On the de Rham-Witt complex attached to a semi-stable family.
Compositio Math. 78 (1991), no. 3, 241--260.

\bibitem{HK}
Hyodo, Osamu; Kato, Kazuya Semi-stable reduction and crystalline cohomology
with logarithmic poles. P\'eriodes $p$-adiques (Bures-sur-Yvette, 1988). 
Ast\'erisque No. 223 (1994), 221--268. 
\bibitem{Il}
Illusie, Luc Complexe de de\thinspace Rham-Witt et cohomologie
cristalline. (French) Ann. Sci. École Norm. Sup. (4) 12 (1979), 
no. 4, 501--661. 

\bibitem{IR} Illusie, Luc; Raynaud, Michel
Les suites spectrales associ\'es au complex de De Rham-Witt.
Publ. Math. IHES No. 57 (1983), 73-212.
\bibitem{KH} Kim, Minhyong; Hain, Richard M.
A De Rham-Witt approach to crystalline rational homotopy theory.
(Preprint 2002). Available as math.AG/0105008.

\bibitem{Mo}Morgan, John W. The rational homotopy theory of smooth, complex projective varieties (following P. Deligne, P. Griffiths, J. Morgan, and D. Sullivan) 
(Invent. Math. 29
  (1975), no. 3, 245--274). S\'eminaire Bourbaki, Vol. 1975/76, 
28\'eme ann\'ee, Exp. No. 475, pp. 69--80. 
Lecture Notes in Math., Vol. 567, Springer, Berlin, 1977. 
\bibitem{NA}
 Navarro Aznar, V. Sur la théorie de Hodge-Deligne. 
 Invent. Math. 90 (1987), no. 1, 11--76. 
\bibitem{Og}
Ogus, Arthur $F$-crystals on schemes with constant log structure. Special issue in honour of Frans Oort. Compositio Math. 97 (1995), no. 1-2, 187--225. 


\bibitem{Wo}
Wojtkowiak, Zdzislaw Cosimplicial objects in algebraic geometry. Algebraic
$K$-theory and algebraic topology (Lake Louise, AB, 1991), 287--327, NATO Adv. Sci. Inst.
Ser. C Math. Phys. Sci., 407, Kluwer Acad. Publ., Dordrecht, 1993.

\end{thebibliography}
\end{document}